\documentstyle[12pt]{article}

\input amssym.def
\input amssym
\newtheorem{definition}{Definition}
\newenvironment{remark}[1]{\vskip 5pt\noindent{\it Remark.\ }#1}{\vskip
5pt}

\newenvironment{proof}[1]{\noindent{\it Proof.\
}#1}{\hskip 3pt $\Box$\vskip 5pt}

\newtheorem{corollary}[definition]{Corollary}
\newtheorem{assertion}[definition]{Assertion}
\newtheorem{conjecture}[definition]{Conjecture}
\newtheorem{proposition}[definition]{Proposition}
\def\limfunc#1{\mbox{\rm #1}\,}
\def\text#1{\mbox{\rm #1}\,}

\begin{document}

\author{{\bf Steven Duplij} \thanks{%
Alexander von Humboldt Fellow} \thanks{%
On leave of absence from {\sl Theory Division, Nuclear Physics Laboratory,
Kharkov State University, KHARKOV 310077, Ukraine}} \thanks{%
E-mail: duplij@physik.uni-kl.de} \\
{\sl Physics Department, University of Kaiserslautern},\\
{\sl Postfach 3049, D-67653 KAISERSLAUTERN},\\
{\sl Germany}}
\title{{\bf NONINVERTIBILITY AND ``SEMI-'' ANALOGS OF (SUPER) MANIFOLDS, FIBER
BUNDLES AND HOMOTOPIES}}
\date{August 19, 1996}
\maketitle

\begin{abstract}
Supersymmetry contains initially noninvertible objects, but it is common to
deal with the invertible ones only, factorizing former in some extent. We
propose to reconsider this ansatz and try to redefine such fundamental
notions as supermanifolds, fiber bundles and homotopies using some weakening
invertibility conditions. The prefix semi- reflects the fact that the
underlying morphisms form corresponding semigroups consisting of a known
group part and a new ideal noninvertible part. We found that the absence of
invertibility gives us the generalization of the cocycle conditions for
transition functions of supermanifolds and fiber bundles in a natural way,
which can lead to construction of noninvertible analogs of \v{C}ech
cocycles. We define semi-homotopies, which can be noninvertible and describe
mappings into the semi-supermanifolds introduced.
\end{abstract}

\newpage

\section{Introduction}

The noninvertible extension of the notion of a supermanifold seems
intuitively natural in connection with the well known theory of Hopf spaces
which in general have no conditions on existence of inverses (in homotopy
sense) \cite{stasheff}. Several guesses were made in the past concerning
inner noninvertibility inherent in the supermanifold theory, e.g. ``...a
general SRS needs not have a body'' \cite{cra/rab}, ``...there may be no
inverse projection (body map \cite{rog1}) at all'' \cite{pen}, or ``...a
body may not even exist in the most extreme examples'' \cite{bry3}. In
particular, while investigating noninvertibility properties of
superconformal symmetry \cite{dup12,dup11} it was assumed \cite{dup6,dup10}
the possible existence of supersymmetric objects analogous to super Riemann
surfaces, but without body, and shown preliminary how to construct them \cite
{dup13}. The noninvertibility in supermanifold theory actually arises from
odd nilpotent elements and zero divisors of underlying Grassmann-Banach
algebras (see \cite{khr1,pes3,she} for nontrivial examples). In the infinite
dimensional case there exist (topologically) quasinilpotent odd elements
which are not really nilpotent \cite{pes6}, and, moreover, in some
superalgebras one can construct pure soul elements which are not nilpotent
even topologically \cite{pes7}. We should also mention the possibility of
definition of a supermanifold without the notion of topological space \cite
{min}. Other supermanifold problems with odd directions (and therefore
connected with noninvertibility in either event) are described in \cite
{bar/bru/her,bry2,cat/rei/teo,khr2,lei2}.

It is well known that the standard patching definition of a supermanifold
\cite{rab/cra2,rog1,vla/vol} does not essentially differ from one of an
ordinary manifold \cite{kosinski,lang}. Mostly the word ``super-'' makes
them different (as however many other definitions too \cite{dewitt}).
Nothing in principle was changed from the time of Gauss who used the concept
of manifold in carthography of earth's surface. And now while constructing
super generalizations of manifolds one thinks about this picture only and
imagine intuitively the earth only, in spite of the fact that in supercase
there are much more abstract possibilities to construct such and similar
objects. The interpretation of the nature of anticommuting variables can be
also dramatically changed in future (see e.g. \cite{yap}). Therefore, it is
logical from the very beginning to weaken (but as fine as possible) the
initial invertibility restrictions. After consistent building of an object
which has a more complicated and nontrivial structure in this general case,
the requirement of invertibility can lead to several nonequivalent
projection including new ones invisible previously.

\section{Standard patch definition of supermanifold}

Let us consider the standard patch definition \cite{rab/cra2,rab1,vla/vol}
of a supermanifold $M$. We cover it by a collection of superdomains $%
U_\alpha $ such that $M=\sum U_\alpha $. Then we take in every domain some
functions ({\sl coordinate maps}) $\varphi _\alpha :U_\alpha \rightarrow
D^{n|m}\subset {\Bbb R}^{n|m}$, where ${\Bbb R}^{n|m}$ is a superspace in
which our super ``ball'' lives and $D^{n|m}$ is an open domain in ${\Bbb R}%
^{n|m}$. Next we call the pair $\left\{ U_\alpha ,\varphi _\alpha \right\} $
a {\sl local chart} and claim that the union of charts $\cup \left\{
U_\alpha ,\varphi _\alpha \right\} $ is an atlas of a {\sl supermanifold}.

Next we introduce gluing {\sl transition functions} as follows. Let $%
U_{\alpha \beta }=U_\alpha \cap U_\beta \neq \emptyset $ and
\begin{equation}  \label{0}
\begin{array}{c}
\varphi _\alpha :U_\alpha \rightarrow V_\alpha \subset {\Bbb R}^{n|m}, \\
\varphi _\beta :U_\beta \rightarrow V_\beta \subset {\Bbb R}^{n|m}.
\end{array}
\end{equation}
Then the above morphisms are restricted to $\varphi _\alpha :U_{\alpha \beta
}\rightarrow V_{\alpha \beta }=V_\alpha \cap \varphi _\alpha \left(
U_{\alpha \beta }\right) $ and $\varphi _\beta :U_{\alpha \beta }\rightarrow
V_{\beta \alpha }=V_\beta \cap \varphi _\beta \left( U_{\alpha \beta
}\right) $. The maps $\Phi_{\alpha \beta }:V_{\beta \alpha }\rightarrow
V_{\alpha \beta }$ which are called to make the following diagram

\begin{equation}  \label{1}
\setlength{\unitlength}{.25in}\begin{picture}(4,4)
\put(0,3){\makebox(1,1){\large ${U_{\alpha\beta}}$}}
\put(3,3){\makebox(1,1){\large ${V_{\beta\alpha}}$}}
\put(3,0){\makebox(1,1){\large ${V_{\alpha\beta}}$}} \put(0.7,2){{\small
$\varphi_\alpha $}} \put(1,3){\vector(1,-1){1.9}} \put(1.6,3.7){{\small
$\varphi_\beta$}} \put(1.25,3.5){\vector(1,0){1.65}}
\put(3.5,2.85){\vector(0,-1){1.75}} \put(3.7,2){{\small
$\Phi_{\alpha\beta}$}} \end{picture}
\end{equation}

\noindent  to be commuted are transition functions of a manifold
in a given atlas. Here we stress, first, that $U_{\alpha \beta }\subset M$
and $V_{\alpha \beta },V_{\beta \alpha }\subset {\Bbb R}^{n|m}$. Second,
from (\ref{1}) one usually concludes that
\begin{equation}  \label{2}
\Phi_{\alpha \beta }=\varphi _\alpha \circ \varphi _\beta ^{-1}
\end{equation}

The transition (super) functions $\Phi _{\alpha \beta }$ give us possibility
to restore the whole (super) manifold from individual charts and coordinate
maps. Indeed they contain all information about the (super) manifold. They
may belong to different functional classes, which gives possibility to
specify more narrow classes of manifolds and supermanifolds, for instance
(super) smooth, analytic, Lipschitz and others \cite{kosinski,okubo}. Mostly
the prefix ``super-'' only distinguishes the patch definitions of a manifold
and supermanifold (which gives us possibility to write it in brackets) and
the properties of $\Phi _{\alpha \beta }$ \cite{bry2,rog1,dewitt}. Here we
do not discuss them in detail and try to put minimum restrictions on $\Phi
_{\alpha \beta }$, concentrating our attention on their abstract properties
and generalizations following from them.

Additionally, from (\ref{2}) it follows that transition functions satisfy
the cocycle conditions

\begin{equation}
\Phi _{\alpha \beta }^{-1}=\Phi _{\beta \alpha }  \label{3}
\end{equation}
on $U_\alpha \cap U_\beta $ and
\begin{equation}
\Phi _{\alpha \beta }\Phi _{\beta \gamma }\Phi _{\gamma \alpha }=1_{\alpha
\alpha }  \label{4}
\end{equation}

\noindent on triple overlaps $U_\alpha \cap U_\beta \cap U_\gamma $, where $%
1_{\alpha \alpha }\stackrel{def}{=}\limfunc{id}\left( U_\alpha \right) $.

Usually it is implied that all the maps $\varphi _\alpha $ are
homeomorphisms, and they can be described by one-to-one invertible
continuous (super) smooth functions (i.e. one wants ``to jump'' in both
directions between any two intersecting domains $U_\alpha \cap U_\beta \neq
\emptyset $). First, it is reasonable (from 19th century carthography
viewpoint) not to distinguish between $U_\alpha $ and $D^{n|m}$, i.e.
locally supermanifolds are as the whole superspace ${\Bbb R}^{n|m}$. For
earth it is right, but for superspace -- questionable. The matter is not
only in more rich fiber bundle \cite{bru/cia1,bru/per2,gra/teo} and sheaf
\cite{kos,leites,lei1} structures due to consideration of all constructions
over underlying Grassmann (or more general \cite{pes7,she,vla/vol,jad/pil})
algebra. The problem lies in another abstract level of the constructions, if
the invertibility conditions are weakened in some extent.

\section{Noninvertible extension of supermanifold}

Earlier there was the following common prescription: one had ready objects
(e.g. real manifolds which can be investigated almost visually), and then
using various methods and guesses one found restrictions on transition
functions. Notwithstanding, noninvertible functions were simply excluded
(saying magic words ``factorizing by nilpotents we again derive the
well-known result'') from consideration, because of desire to be in the
nearest analogy with intuitively clear and understandable nonsupersymmetric
case.

Here we go in opposite direction: we know that in supermathematics
noninvertible variables and functions {\it do exist}. Which objects could be
constructed by means of them? What gives ``factorizing by non-nilpotents'',
i.e. consideration of non-group features of theory? How changes the general
abstract meaning of the most important notions, e.g. manifolds and fiber
bundles? We now try to leave aside inner structure of noninvertible objects
analogous to supermanifolds and concentrate our attention on there general
abstract definitions. Further we think to work them out in more detail.

\subsection{Problem of division in superanalysis}

We begin with old division problem of superanalysis \cite{berezin,rab/cra1}.
Can we solve the simplest equation $ax=b$ with even noninvertible $a$ , i.e.
when $a^n=0$ ? Here $a,b,x\in \Lambda _0$ , where $\Lambda $ is a
commutative Banach ${\Bbb Z}_2$-graded superalgebra over a field ${\Bbb K}$
(where ${\Bbb K}={\Bbb R}{\Bbb ,}$ or ${\Bbb C}$) with a decomposition into
the direct sum: $\Lambda =\Lambda _0\oplus \Lambda _1$ (the elements $a$
from $\Lambda _0$ and $\Lambda _1$ are homogeneous and have the fixed even
and odd parity defined as $\left| a\right| \stackrel{def}{=}\left\{ i\in
\left\{ 0,1\right\} ={\Bbb Z}_2|\,a\in \Lambda _i\right\} $ ). The answer is
``yes'': we should expand both sides in series of generators of $\Lambda $
and equal the coefficients before the same power of generators. The simplest
example is $\alpha x=2\alpha \beta \gamma $ which has the solution $x=2\beta
\gamma +\limfunc{Ann}\alpha $. This means that the extended division over
noninvertible variable {\it can be defined}, i.e. the solution of $ax=b$
when $a^n=0$ {\it exists}, but it can be found by indirect way (without
notorious factorization by nilpotents). But how to solve $\alpha x=2\beta $
or $\alpha \beta x=2\gamma \rho $ ? For instance, one can use the technique
of the nilpotent semigroup theory \cite{gri1,lal,sul1}. Another possibility
to solve such equations is to exploit the abstract semigroup methods \cite
{ljapin} and define solutions as equivalent classes of variables or
functions. We only mark now that the solution of $ax=b$ exists when $b$ is
also nilpotent $b^k=0$. The same conclusion can be applied to functional
equations containing superfunctions: they {\it can} be solved, and so they
{\it have to be considered on a par} with invertible ones. We stress that
among ordinary functions there exist noninvertible ones as well \cite
{mag11,mag6}, but the kind of noninvertibility considered here is very
special: it appears only due to the existence of nilpotents in underlying
superalgebra \cite{khr1,pes6,she}. Here we do not consider concrete
equations and methods of their solving, we only use the fact of their
existence to reformulate some definitions and extend well-known notions.

\subsection{Definition of a semi-supermanifold}

Now we formulate a patch definition of an object analogous to supermanifold,
i.e. try to weaken demand of invertibility of coordinate maps (\ref{0}). Let
us consider a (super) space $M$ covered by open sets $U_\alpha $ as $M=\sum
U_\alpha $ . Let the maps $\varphi _\alpha :U_\alpha \rightarrow V_\alpha
\subset {\Bbb R}^{n|m}$ are not all homeomorphisms, i.e. among them there
are noninvertible maps.

\begin{definition}
A{\sl \ chart }is a pair $\left( U_{\alpha ,}\varphi _\alpha \right) $
,where $\varphi _\alpha $ is invertible. A{\sl \ semi-chart }is a pair $%
\left( \tilde{U}_{\alpha ,}\tilde{\varphi}_\alpha \right) $ $,$ where $%
\tilde{\varphi}_\alpha $ is noninvertible.
\end{definition}

\begin{definition}
A{\sl \ semi-atlas }is a union of charts and semi-charts.
\end{definition}

\begin{definition}
A{\sl \ semi-supermanifold }is a superspace $M$ represented as a semi-atlas.
\end{definition}

How to define an analog of transition functions? We should consider the same
diagram (\ref{1}), but we cannot use (\ref{2}) through noninvertibility of
some of $\varphi _\alpha $'s.

\begin{definition}
Gluing {\sl semi-transition functions }of a semi-supermanifold are defined
by the equations
\begin{equation}
\Phi _{\alpha \beta }\circ \varphi _\beta =\varphi _\alpha   \label{5}
\end{equation}
\begin{equation}
\Phi _{\beta \alpha }\circ \varphi _\alpha =\varphi _\beta   \label{6}
\end{equation}
\end{definition}

We stress that to determine $\Phi_{\alpha \beta }$ the equation (\ref{5})
cannot be solved by (\ref{2}). Instead we should find artificial methods of
its solving, e.g. as in previous subsection, expanding in superalgebra
generator series, or using abstract semigroup methods \cite{ljapin} which
consider the solutions of noninvertible equations as equivalence classes.

The functions $\Phi _{\beta \alpha }$ are now determined not from (\ref{3})
in which the left hand side is not well defined, but from the commutative
diagram

\begin{equation}
\setlength{\unitlength}{.25in}\begin{picture}(4,4)
\put(0,3){\makebox(1,1){\large ${U_{\alpha\beta}}$}}
\put(3,3){\makebox(1,1){\large ${V_{\beta\alpha}}$}}
\put(3,0){\makebox(1,1){\large ${V_{\alpha\beta}}$}} \put(0.7,2){{\small
$\varphi_\alpha $}} \put(1,3){\vector(1,-1){1.9}} \put(1.6,3.7){{\small
$\varphi_\beta$}} \put(1.15,3.5){\vector(1,0){1.7}}
\put(3.5,1.1){\vector(0,1){1.75}} \put(3.7,2){{\small $\Phi_{\beta\alpha}$}}
\end{picture}  \label{7}
\end{equation}

\noindent and the equation (\ref{6}) following from it. They are also can be
noninvertible and therefore the cocycle conditions should be modified not to
use only invertible functions.

\begin{remark}
Even in the standard case the cocycle conditions (\ref{4}) for
supermanifolds are not automatically satisfied when (\ref{2}) holds, and
therefore they should be imposed by hand \cite{nel}.
\end{remark}

So instead of (\ref{3}) and (\ref{4}) we have

\begin{conjecture}
The semi-transition functions of a semi-supermanifold satisfy the following
relations
\begin{equation}
\Phi _{\alpha \beta }\circ \Phi _{\beta \alpha }\circ \Phi _{\alpha \beta
}=\Phi _{\alpha \beta }  \label{8}
\end{equation}

\noindent on $U_\alpha \cap U_\beta $ overlaps and
\begin{equation}
\Phi _{\alpha \beta }\circ \Phi _{\beta \gamma }\circ \Phi _{\gamma \alpha
}\circ \Phi _{\alpha \beta }=\Phi _{\alpha \beta },  \label{90}
\end{equation}
\begin{equation}
\Phi _{\beta \gamma }\circ \Phi _{\gamma \alpha }\circ \Phi _{\alpha \beta
}\circ \Phi _{\beta \gamma }=\Phi _{\beta \gamma },  \label{91}
\end{equation}
\begin{equation}
\Phi _{\gamma \alpha }\circ \Phi _{\alpha \beta }\circ \Phi _{\beta \gamma
}\circ \Phi _{\gamma \alpha }=\Phi _{\gamma \alpha }  \label{92}
\end{equation}
on triple overlaps $U_\alpha \cap U_\beta \cap U_\gamma $ and
\begin{equation}
\Phi _{\alpha \beta }\circ \Phi _{\beta \gamma }\circ \Phi _{\gamma \rho
}\circ \Phi _{\rho \alpha }\circ \Phi _{\alpha \beta }=\Phi _{\alpha \beta },
\label{90w}
\end{equation}
\begin{equation}
\Phi _{\beta \gamma }\circ \Phi _{\gamma \rho }\circ \Phi _{\rho \alpha
}\circ \Phi _{\alpha \beta }\circ \Phi _{\beta \gamma }=\Phi _{\beta \gamma
},  \label{91w}
\end{equation}
\begin{equation}
\Phi _{\gamma \rho }\circ \Phi _{\rho \alpha }\circ \Phi _{\alpha \beta
}\circ \Phi _{\beta \gamma }\circ \Phi _{\gamma \rho }=\Phi _{\gamma \rho },
\label{92v}
\end{equation}
\begin{equation}
\Phi _{\rho \alpha }\circ \Phi _{\alpha \beta }\circ \Phi _{\beta \gamma
}\circ \Phi _{\gamma \rho }\circ \Phi _{\rho \alpha }=\Phi _{\rho \alpha }
\label{92w}
\end{equation}
\noindent on $U_\alpha \cap U_\beta \cap U_\gamma \cap U_\rho $ .
\end{conjecture}

Here the first relation (\ref{8}) is called to generalize the first cocycle
condition (\ref{3}), while other relations correspond (\ref{4}). We call (%
\ref{8})--(\ref{92w}) {\sl tower relations}.

\begin{definition}
A semi-supermanifold is called {\sl reflexive} if, in addition to (\ref{8}%
)--(\ref{92w}), the semi-transition functions satisfy to the reflexivity
conditions
\begin{equation}
\Phi _{\beta \alpha }\circ \Phi _{\alpha \beta }\circ \Phi _{\beta \alpha
}=\Phi _{\beta \alpha }  \label{8a}
\end{equation}

\noindent on $U_\alpha \cap U_\beta $ overlaps and
\begin{equation}
\Phi _{\alpha \gamma }\circ \Phi _{\gamma \beta }\circ \Phi _{\beta \alpha
}\circ \Phi _{\alpha \gamma }=\Phi _{\alpha \gamma },  \label{9a}
\end{equation}
\begin{equation}
\Phi _{\gamma \beta }\circ \Phi _{\beta \alpha }\circ \Phi _{\alpha \gamma
}\circ \Phi _{\gamma \beta }=\Phi _{\gamma \beta },  \label{9a1}
\end{equation}
\begin{equation}
\Phi _{\beta \alpha }\circ \Phi _{\alpha \gamma }\circ \Phi _{\gamma \beta
}\circ \Phi _{\beta \alpha }=\Phi _{\beta \alpha }  \label{9a2}
\end{equation}
on triple overlaps $U_\alpha \cap U_\beta \cap U_\gamma $ and
\begin{equation}
\Phi _{\alpha \rho }\circ \Phi _{\rho \gamma }\circ \Phi _{\gamma \beta
}\circ \Phi _{\beta \alpha }\circ \Phi _{\alpha \rho }=\Phi _{\alpha \rho },
\label{9aw}
\end{equation}
\begin{equation}
\Phi _{\rho \gamma }\circ \Phi _{\gamma \beta }\circ \Phi _{\beta \alpha
}\circ \Phi _{\alpha \rho }\circ \Phi _{\rho \gamma }=\Phi _{\rho \gamma },
\label{9aw1}
\end{equation}
\begin{equation}
\Phi _{\gamma \beta }\circ \Phi _{\beta \alpha }\circ \Phi _{\alpha \rho
}\circ \Phi _{\rho \gamma }\circ \Phi _{\gamma \beta }=\Phi _{\gamma \beta },
\label{9a1w}
\end{equation}
\begin{equation}
\Phi _{\beta \alpha }\circ \Phi _{\alpha \rho }\circ \Phi _{\rho \gamma
}\circ \Phi _{\gamma \beta }\circ \Phi _{\beta \alpha }=\Phi _{\beta \alpha }
\label{9a2w}
\end{equation}
\noindent on $U_\alpha \cap U_\beta \cap U_\gamma \cap U_\rho $ .
\end{definition}

\begin{remark}
One could think that the reflexivity conditions differ from (\ref{8})--(\ref
{92w}) by index permutations only. This is true. But the functions $\Phi
_{\alpha \beta }$ entering in (\ref{8})--(\ref{92w}) and in (\ref{8a})--(\ref
{9a2w}) {\it are the same}, therefore the latter are independent equations
imposed on $\Phi _{\alpha \beta }$.
\end{remark}

\begin{assertion}
The relations analogous to (\ref{8})--(\ref{9a2w}), but having two or more
multipliers in the right hand side are consequences of them.
\end{assertion}

\begin{proof}
For instance, consider
\begin{equation}
\Phi _{\alpha \beta }\circ \Phi _{\beta \gamma }\circ \Phi _{\gamma \alpha
}\circ \Phi _{\alpha \beta }\circ \Phi _{\beta \gamma }=\Phi _{\alpha \beta
}\circ \Phi _{\beta \gamma }.  \label{10}
\end{equation}
Multiplying from the right on $\Phi _{\alpha \beta }$ we derive
\begin{equation}
\Phi _{\alpha \beta }\circ \Phi _{\beta \gamma }\circ \Phi _{\gamma \alpha
}\circ \Phi _{\alpha \beta }\circ \Phi _{\beta \gamma }\circ \Phi _{\alpha
\beta }=\Phi _{\alpha \beta }\circ \Phi _{\beta \gamma }\circ \Phi _{\alpha
\beta }.  \label{10p}
\end{equation}
Then using (\ref{8}) we obtain
\begin{equation}
\Phi _{\alpha \beta }\circ \Phi _{\beta \gamma }\circ \Phi _{\gamma \alpha
}\circ \Phi _{\alpha \beta }=\Phi _{\alpha \beta }  \label{10q}
\end{equation}
which coincides with (\ref{90}).
\end{proof}

\begin{remark}
In any actions with noninvertible functions $\Phi _{\alpha \beta }$ we are
not allowed to cancellate, because the semigroup of $\Phi _{\alpha \beta }$%
's is a semigroup without cancellation, and we are forced to exploit the
corresponding methods \cite{csa/thu,hof,mag4,vec}.
\end{remark}

\begin{corollary}
The relations (\ref{8})--(\ref{9a2w}) satisfy identically in the standard
invertible case, i.e. when the conditions (\ref{2}),(\ref{3}) and (\ref{4})
hold valid.
\end{corollary}

So that, $\Phi _{\alpha \beta }$ satisfying the relations (\ref{8})--(\ref
{9a2w}) can be viewed as some noninvertible generalization of the transition
functions as cocycles in the \v {C}ech cohomology of coverings \cite
{maclane,switzer}.

\subsection{Orientation of semi-manifolds}

It is well known that orientation of ordinary manifolds is determined by the
Jacobian sign of transition functions $\Phi _{\alpha \beta }$ written in
terms of local coordinates on $U_\alpha \cap U_\beta $ overlaps \cite
{kosinski,schwarz}. Since this sign belong to ${\Bbb Z}_2$ , there exist two
orientations on $U_\alpha $. Two overlapping charts are {\sl consistently
oriented} (or {\sl orientation preserving}) if $\Phi _{\alpha \beta }$ has
positive Jacobian, and a manifold is {\sl orientable} if it can be covered
by such charts, thus there are two kinds of manifolds: orientable and
nonorientable. In supersymmetric case the role of Jacobian plays Berezinian
\cite{berezin} which has a ``sign'' belonging to ${\Bbb Z}_2\oplus {\Bbb Z}%
_2 $ \cite{voronov,tuy}, and so there are four orientations on $U_\alpha $
and five corresponding kinds of supermanifold orientability \cite{sha1}.

\begin{definition}
In case a nonvanishing Berezinian of $\Phi _{\alpha \beta }$ is nilpotent
(and so has no definite sign in the previous sense) there exists additional
{\sl nilpotent orientation }on $U_\alpha $ of a semi-supermanifold.
\end{definition}

A degree of nilpotency of Berezinian allows us to classify
semi-supermanifolds having nilpotent orientability.

\subsection{Obstructedness and semi-manifolds}

The semi-supermanifolds defined above belong to a class of so called
obstructed semi-manifolds. Let us rewrite (\ref{2}), (\ref{3}) and (\ref{4})
as the following (infinite) series
\begin{equation}
\begin{array}{cc}
n=1: & \Phi _{\alpha \alpha }=1_{\alpha \alpha },
\end{array}
\label{11a}
\end{equation}
\begin{equation}
\begin{array}{cc}
n=2: & \Phi _{\alpha \beta }\circ \Phi _{\beta \alpha }=1_{\alpha \alpha },
\end{array}
\label{11b}
\end{equation}
\begin{equation}
\begin{array}{cc}
n=3: & \Phi _{\alpha \beta }\circ \Phi _{\beta \gamma }\circ \Phi _{\gamma
\alpha }=1_{\alpha \alpha },
\end{array}
\label{11c}
\end{equation}
\begin{equation}
\begin{array}{cc}
n=4: & \Phi _{\alpha \beta }\circ \Phi _{\beta \gamma }\circ \Phi _{\gamma
\delta }\circ \Phi _{\delta \alpha }=1_{\alpha \alpha }
\end{array}
\label{11d}
\end{equation}
\[
\begin{array}{cc}
\cdots  & \cdots
\end{array}
\]

\begin{definition}
A semi-manifold is called {\sl obstructed }if some of the cocycle conditions
(\ref{11a})--(\ref{11d}) are broken.
\end{definition}

\begin{remark}
The introduced notion of obstructed manifold should not be mixed with the
notion of obstruction for ordinary manifolds \cite{baues} and supermanifolds
\cite{berezin} or obstruction to extensions \cite{maclane} and in the theory
of characteristic classes \cite{mil/sta,kam/ton}.
\end{remark}

It can happen that starting from some $n=n_m$ all higher cocycle conditions
hold valid.

\begin{definition}
{\sl Obstructedness degree} of a semi-manifold is a maximal $n_m$ for which
the cocycle conditions (\ref{11a})--(\ref{11d}) are broken. If all of them
hold valid, then $n_m\stackrel{def}{=}0$.
\end{definition}

\begin{corollary}
Ordinary manifolds (with invertible transition functions) have vanishing
obstructedness, and the obstructedness degree for them is equal to zero,
i.e. $n_m=0$.
\end{corollary}

\begin{conjecture}
The obstructed semi-manifolds may have non-vanishing ordinary obstruction
which can be calculated extending the standard methods \cite{berezin} to the
non-invertible case.
\end{conjecture}

Therefore, using the obstructedness degree $n_m$, we have possibility to
classify semi-manifolds properly.

\subsection{Tower identity semigroup}

Let us consider a series of the selfmaps ${\rm e}_{\alpha \alpha }^{\left(
n\right) }:U_\alpha \rightarrow U_\alpha $ of a semi-manifold defined as
\begin{equation}
{\rm e}_{\alpha \alpha }^{\left( 1\right) }=\Phi _{\alpha \alpha },
\label{e11a}
\end{equation}
\begin{equation}
{\rm e}_{\alpha \alpha }^{\left( 2\right) }=\Phi _{\alpha \beta }\circ \Phi
_{\beta \alpha },  \label{e11b}
\end{equation}
\begin{equation}
{\rm e}_{\alpha \alpha }^{\left( 3\right) }=\Phi _{\alpha \beta }\circ \Phi
_{\beta \gamma }\circ \Phi _{\gamma \alpha },  \label{e11c}
\end{equation}
\begin{equation}
{\rm e}_{\alpha \alpha }^{\left( 4\right) }=\Phi _{\alpha \beta }\circ \Phi
_{\beta \gamma }\circ \Phi _{\gamma \delta }\circ \Phi _{\delta \alpha }
\label{e11d}
\end{equation}
\[
\begin{array}{cc}
\cdots & \cdots
\end{array}
\]

We will call ${\rm e}_{\alpha \alpha }^{\left( n\right) }$'s {\sl tower
identities}. From (\ref{11a})--(\ref{11d}) it follows

\begin{assertion}
For an ordinary supermanifolds all tower identities coincide with the usual
identity map
\begin{equation}
{\rm e}_{\alpha \alpha }^{\left( n\right) }=1_{\alpha \alpha }.  \label{e1}
\end{equation}
\end{assertion}

The obstructedness degree can be treated as a maximal $n=n_m$ for which
tower identities differ from the identity, i.e. (\ref{e1}) is broken. So the
tower identities give the measure of distinction of a semi-supermanifold
from an ordinary supermanifold. Being an important inner characteristic the
tower identities play a deep fundamental role in description of
semi-supermanifolds. Therefore, we will study some of their properties in
detail.

\begin{proposition}
The tower identities are units for the semi-transition functions
\begin{equation}
{\rm e}_{\alpha \alpha }^{\left( n\right) }\circ \Phi _{\alpha \beta }=\Phi
_{\alpha \beta },  \label{ef1}
\end{equation}
\begin{equation}
\Phi _{\alpha \beta }\circ {\rm e}_{\beta \beta }^{\left( n\right) }=\Phi
_{\alpha \beta }.  \label{ef2}
\end{equation}
\end{proposition}

\begin{proof}
It follows directly from the tower relations (\ref{8})--(\ref{92w}) and the
definition (\ref{8})--(\ref{92w}).
\end{proof}

\begin{proposition}
The tower identities are idempotents
\begin{equation}
{\rm e}_{\alpha \alpha }^{\left( n\right) }\circ {\rm e}_{\alpha \alpha
}^{\left( n\right) }={\rm e}_{\alpha \alpha }^{\left( n\right) }.
\label{eee}
\end{equation}
\end{proposition}

\begin{proof}
We prove the statement for $n=2$ and for other $n$ it can be proved by
induction. We write (\ref{eee}) as
\[
{\rm e}_{\alpha \alpha }^{\left( 2\right) }\circ {\rm e}_{\alpha \alpha
}^{\left( 2\right) }={\rm e}_{\alpha \alpha }^{\left( 2\right) }\circ \Phi
_{\alpha \beta }\circ \Phi _{\beta \alpha }=\left( {\rm e}_{\alpha \alpha
}^{\left( 2\right) }\circ \Phi _{\alpha \beta }\right) \circ \Phi _{\beta
\alpha }.
\]

Then using (\ref{ef1}) we obtain
\[
\left( {\rm e}_{\alpha \alpha }^{\left( 2\right) }\circ \Phi _{\alpha \beta
}\right) \circ \Phi _{\beta \alpha }=\Phi _{\alpha \beta }\circ \Phi _{\beta
\alpha }={\rm e}_{\alpha \alpha }^{\left( 2\right) }.
\]
\end{proof}

The functional nonsupersymmetric equations of the above
kind were studied in \cite{bed/wal}.

\begin{definition}
{\sl Conjugate tower identities }correspond to the same partition of the
semi-supermanifold and consists of the semi-transition functions taken in
opposite order
\begin{equation}
{\rm \tilde{e}}_{\alpha \alpha }^{\left( 1\right) }={\rm e}_{\alpha \alpha
}^{\left( 1\right) },  \label{ee11a}
\end{equation}
\begin{equation}
{\rm \tilde{e}}_{\alpha \alpha }^{\left( 2\right) }={\rm e}_{\alpha \alpha
}^{\left( 2\right) },  \label{ee11b}
\end{equation}
\begin{equation}
{\rm \tilde{e}}_{\alpha \alpha }^{\left( 3\right) }=\Phi _{\alpha \gamma
}\circ \Phi _{\gamma \beta }\circ \Phi _{\beta \alpha },  \label{ee11c}
\end{equation}
\begin{equation}
{\rm \tilde{e}}_{\alpha \alpha }^{\left( 4\right) }=\Phi _{\alpha \delta
}\circ \Phi _{\delta \gamma }\circ \Phi _{\gamma \beta }\circ \Phi _{\beta
\alpha }  \label{ee11d}
\end{equation}
\[
\begin{array}{cc}
\cdots  & \cdots
\end{array}
\]
\end{definition}

The conjugate tower identities play the role of tower identities, but for
reflexivity conditions (\ref{8a})--(\ref{9a2w}). By analogy with (\ref{ef1}%
)--(\ref{ef2}) we have

\begin{proposition}
The conjugate tower identities are {\sl reflexive units} for the
semi-transition functions
\begin{equation}
{\rm \tilde{e}}_{\beta \beta }^{\left( n\right) }\circ \Phi _{\beta \alpha
}=\Phi _{\beta \alpha },  \label{eef1}
\end{equation}
\begin{equation}
\Phi _{\beta \alpha }\circ {\rm \tilde{e}}_{\alpha \alpha }^{\left( n\right)
}=\Phi _{\beta \alpha }.  \label{eef2}
\end{equation}
\end{proposition}

\begin{assertion}
For the same partition the conjugate tower identities annihilate the tower
identities in the following sense
\begin{equation}
{\rm e}_{\alpha \alpha }^{\left( n\right) }\circ {\rm \tilde{e}}_{\alpha
\alpha }^{\left( n\right) }={\rm e}_{\alpha \alpha }^{\left( 2\right) }.
\label{eeea}
\end{equation}
\end{assertion}

\begin{proof}
Let us consider the case $n=3$. Using the definitions we derive
\[
{\rm e}_{\alpha \alpha }^{\left( 3\right) }\circ {\rm \tilde{e}}_{\alpha
\alpha }^{\left( 3\right) }=\Phi _{\alpha \beta }\circ \Phi _{\beta \gamma
}\circ \Phi _{\gamma \alpha }\circ \Phi _{\alpha \gamma }\circ \Phi _{\gamma
\beta }\circ \Phi _{\beta \alpha }
\]
\[
=\Phi _{\alpha \beta }\circ \Phi _{\beta \gamma }\circ \left( \Phi _{\gamma
\alpha }\circ \Phi _{\alpha \gamma }\right) \circ \Phi _{\gamma \beta }\circ
\Phi _{\beta \alpha }=\Phi _{\alpha \beta }\circ \Phi _{\beta \gamma }\circ
{\rm e}_{\gamma \gamma }^{\left( 2\right) }\circ \Phi _{\gamma \beta }\circ
\Phi _{\beta \alpha }
\]
\[
=\Phi _{\alpha \beta }\circ \left( \Phi _{\beta \gamma }\circ \Phi _{\gamma
\beta }\right) \circ \Phi _{\beta \alpha }=\Phi _{\alpha \beta }\circ {\rm e}%
_{\beta \beta }^{\left( 2\right) }\circ \Phi _{\beta \alpha }=\Phi _{\alpha
\beta }\circ \Phi _{\beta \alpha }={\rm e}_{\alpha \alpha }^{\left( 2\right)
}.
\]

For other $n$ the statement can be proved by induction.
\end{proof}

\begin{definition}
A semi-supermanifold is {\sl nice} if the tower identities do not depend on
a given partition.
\end{definition}

The multiplication of the tower identities of a nice semi-supermanifold can
be defined as follows
\begin{equation}
{\rm e}_{\alpha \alpha }^{\left( n\right) }\circ {\rm e}_{\alpha \alpha
}^{\left( m\right) }={\rm e}_{\alpha \alpha }^{\left( n+m\right) }.
\label{eee1}
\end{equation}

\begin{assertion}
The multiplication (\ref{eee1}) is associative.
\end{assertion}

Therefore, we are able to give

\begin{definition}
The tower identities of a nice semi-supermanifold form {\sl a tower semigroup%
} under the multiplication (\ref{eee1}).
\end{definition}

So we obtained the quantitative description of inner noninvertibility
properties of semi-supermanifolds. The introduced tower semigroup plays the
same role for semi-supermanifolds as the fundamental group for ordinary
manifolds \cite{fom/fuc/gut,maclane,switzer}.

\subsection{Semicommutative diagrams and $n$-regularity}

The above constructions have the general importance for any set of
noninvertible mappings.

The extension of $n=2$ cocycle given by (\ref{8}) can be viewed as some
analogy with regular elements in semigroups \cite{cli/pre1,petrich3} or
generalized inverses in matrix theory \cite{pen1,rao/mit}, category theory
\cite{dav/rob} and theory of generalized inverses of morphisms \cite{nashed}%
.

The relations (\ref{90})--(\ref{92w}) and with other $n$ can be considered
as noninvertible analogue of regularity for higher cocycles. Therefore, by
analogy with (\ref{8})--(\ref{92w}) it is natural to formulate the general

\begin{definition}
An noninvertible mapping $\Phi _{\alpha \beta }$ is $n${\sl -regular}, if it
satisfies the following conditions
\begin{equation}
\stackrel{n}{\overbrace{\Phi _{\alpha \beta }\circ \Phi _{\beta \gamma
}\circ \ldots \circ \Phi _{\rho \alpha }\circ \Phi _{\alpha \beta }}}=\Phi
_{\alpha \beta }  \label{88}
\end{equation}
\[
+\, perm
\]

on overlaps $\stackrel{n}{\overbrace{U_\alpha \cap U_\beta \cap \ldots \cap
U_\rho }}$.
\end{definition}

The formula (\ref{8}) describes $3$-regular mappings, the relations (\ref{90}%
)--(\ref{92}) correspond to $4$-regular ones, and (\ref{90w})--(\ref{92w})
give $5$-regular mappings.

\begin{remark}
The $3$-regularity coincides with the ordinary regularity.
\end{remark}

Another definition of $n$-regularity can be given by the formulas (\ref{ef1}%
)--(\ref{ef2}).

The higher regularity conditions change dramatically the general diagram
technique of morphisms, when we turn to noninvertible ones. Indeed, the
commutativity of invertible morphism diagrams is based on the relations (\ref
{11a})--(\ref{11d}), i.e. on the fact that the tower identities are ordinary
identities (\ref{e1}). When morphisms are noninvertible (a
semi-supermanifold has a nonvanishing obstructedness), we cannot ``return to
the same point'', because ${\rm e}_{\alpha \alpha }^{\left( n\right) }\neq
1_{\alpha \alpha }$, and we have to consider ``nonclosed'' diagrams due to
the fact that the relation ${\rm e}_{\alpha \alpha }^{\left( n\right) }\circ
\Phi _{\alpha \beta }=\Phi _{\alpha \beta }$ is noncancellative now.

Summarizing the above statements we propose the following intuitively
consistent changing of the standard diagram technique as applied to
noninvertible morphisms. In every case we get a new arrow which corresponds
to the additional multiplier in   (\ref{ef1}).
Thus, for $n=2$ we obtain

\setlength{\unitlength}{.3in}
\begin{picture}(15,3.5)
\put(4.5,2.8){\makebox(1,1){\bf Invertible}}
\put(4.6,2.2){{ ${\Phi }_{\alpha \beta}$}}
\put(4.6,0.7){{ ${\Phi }_{\beta \alpha }$}}
\put(4,1.8){\vector(1,0){2.1}}
\put(6.1,1.4){\vector(-1,0){2.1}}
\put(8,1.5){$\Longrightarrow $}
\put(10.5,2.8){\makebox(1,1){\bf Noninvertible}}
\put(10.6,0.6){{ ${\Phi }_{\beta \alpha }$}}
\put(10.6,2.2){{ ${\Phi }_{\alpha \beta}$}}
\put(10,1.9){\vector(1,0){2.1}}
\put(12.1,1.3){\vector(-1,0){2.1}}
\put(10,1.6){\vector(1,0){2.1}}
\put(0,1.5){\large \it n=2}
\end{picture}

\noindent which describes the transition from  (\ref{11b})
to   (\ref{8}) and
presents the ordinary regularity condition for morphisms
\cite{dav/rob,nashed}.
The most
intriguing semicommutative diagram is the triangle one

\setlength{\unitlength}{.3in}
\begin{picture}(15,5)
\put(5.8,1.5){\vector(-1,1){1.68}}
\put(4.6,4.2){{ ${\Phi }_{\alpha \beta}$}}
\put(3.5,2){{ ${\Phi}_{\gamma \alpha }$}}
\put(4,3.5){\vector(1,0){2.1}}
\put(6.1,3.2){\vector(0,-1){2}}
\put(8,2.5){$\Longrightarrow $}
\put(14,2.5){+ perm}
\put(6.5,2){{ ${\Phi }_{\beta \gamma}$}}
\put(9.5,2){{ ${\Phi}_{\gamma \alpha }$}}
\put(11.8,1.5){\vector(-1,1){1.68}}
\put(10.6,4.2){{ ${\Phi }_{\alpha \beta}$}}
\put(10,3.8){\vector(1,0){2.1}}
\put(10,3.4){\vector(1,0){2.1}}
\put(12.1,3.2){\vector(0,-1){2}}
\put(12.5,2){{ ${\Phi}_{\beta \gamma}$}}
\put(0,2.5){\large \it n=3}
\end{picture}

\noindent which generalizes the cocycle condition   (\ref{4}).
By analogy one can
write higher order diagrams.

\section{Noninvertibility and semi-bundles}

A similar approach can be applied to the noninvertible extension of fiber
bundles, while defining them globally in terms of open coverings and
transition functions.

Following the standard definitions \cite{mil/sta,okubo,schwarz} and
weakening invertibility we now construct new objects analogous to fiber
bundles.

\subsection{Definition of semi-bundles}

Let $E$ and $M$ be a {\sl total} ({\sl bundle}) superspace and {\sl base}
semi-supermanifold respectively, and $\pi :E\rightarrow M$ be a {\sl %
semi-projection map} which is not necessarily invertible (but can be
smooth). Denote by $F_b$ the set of points of $E$ that are mapped to $b\in M$
(a {\sl pre-image }of $b$ ), i.e. the {\sl semi-fiber over }$b$ is $F_b%
\stackrel{def}{=}\left\{ x\in E\,|\,\pi \left( x\right) =b\right\} $. Then, $%
F=\cup F_b$ is a {\sl semi-fiber}.

\begin{definition}
A {\sl semi-bundle} is ${\cal L}\stackrel{def}{=}\left( E,M,F,\pi \right) $ .
\end{definition}

A {\sl section} $s:M\rightarrow F$ of the fiber bundle $\left( E,M,F,\pi
\right) $ is usually defined by $\pi \left( s\left( b\right) \right) =b$
which in the form $\pi \circ s=1_M$ is very similar to (\ref{3}) and (\ref
{11b}) and holds valid for invertible maps $\pi $ and $s$ only. Therefore, a
very few ordinary nontrivial fiber bundles admit corresponding sections \cite
{mil/sta}.

Thus, using analogy with (\ref{8}), we come to the following

\begin{definition}
A {\sl semi-section} of the semi-bundle ${\cal L}=\left( E,M,F,\pi \right) $
is defined by
\begin{equation}
\pi \circ s\circ \pi =\pi .  \label{12a}
\end{equation}

A {\sl reflexive} semi-section satisfies to the additional condition
\begin{equation}
s\circ \pi \circ s=s.  \label{12b}
\end{equation}
\end{definition}

Let $\tilde{\pi}:M\times F\rightarrow M$ is the canonical semi-projection on
the first factor $\tilde{\pi}\left( b,f\right) =b,\;\,f\in F$, then $\tilde{%
\pi}$ gives rise to a {\sl product fiber bundle}. If $\lambda :E\rightarrow
M\times F$ is a morphism (called a {\sl trivialization}), then $\tilde{\pi}%
\circ \lambda =\pi $ , and the semi-bundle ${\cal L}=\left( E,M,F,\pi
\right) $ is {\sl trivial}. If there exists a continuous map $\eta
:M\rightarrow F$ , then the semi-bundle $\left( M\times F,M,F,\tilde{\pi}%
\right) $ admits the section $s:M\rightarrow M\times F$ given by $s\left(
b\right) =\left( b,\eta \left( b\right) \right) $.

Let $E_\alpha \stackrel{def}{=}\left\{ x\in E\,|\,\pi _\alpha \left(
x\right) =b,\,b\in U_\alpha \subset M\right\} $ (here we do not use the
standard notion $\pi ^{-1}\left( U_\alpha \right) $ for $E_\alpha $
intentionally, because now $\pi _\alpha $ is allowed to be noninvertible),
where $\pi _\alpha :E_\alpha \rightarrow U_\alpha $ is a restriction, i.e. $%
\pi _\alpha \stackrel{def}{=}\pi \mid _{U_\alpha }$ . Then the semi-bundle $%
{\cal L}=\left( E,M,F,\pi \right) $, is {\sl locally trivial}, if $\forall
b\in M\,\exists U_\alpha \ni b$ such that there exists the trivializing
morphisms $\lambda _\alpha :E_\alpha \rightarrow U_\alpha \times F$
satisfying $\tilde{\pi}_\alpha \circ \lambda _\alpha =\pi _\alpha $. That
is, the diagram

\begin{equation}
\setlength{\unitlength}{.25in}\begin{picture}(4,4)
\put(0,3){\makebox(1,1){\large ${E_{\alpha}}$}}
\put(3.8,3){\makebox(1,1){\large ${U_{\alpha} \times F}$}}
\put(3,0){\makebox(1,1){\large ${U_{\alpha}}$}} \put(0.7,2){{\small
$\pi_\alpha $}} \put(1,3){\vector(1,-1){1.9}} \put(1.6,3.7){{\small
$\lambda_\alpha $}} \put(1.1,3.5){\vector(1,0){1.7}}
\put(3.5,2.85){\vector(0,-1){1.75}} \put(3.7,2){{\small $\tilde{\pi}_\alpha
$}} \end{picture}  \label{p3}
\end{equation}

\noindent  commutes.

\begin{definition}
A {\sl semi-section} of a locally trivial semi-bundle ${\cal L}$ is given by
the maps $s_\alpha :U_\alpha \rightarrow E$ which satisfy the compatibility
conditions
\begin{equation}
\lambda _\alpha \circ s_\alpha \mid _b=\lambda _\beta \circ s_\beta \mid
_b,\;\,b\in U_\alpha \cap U_\beta .  \label{13}
\end{equation}
\end{definition}

Now let $\left\{ U_\alpha ,\lambda _\alpha \right\} $ be a trivializing
covering of $\pi $ such that $\cup U_\alpha =M$ and $U_\alpha \cap U_\beta
\neq \emptyset \Rightarrow E_\alpha \cap E_\beta \neq \emptyset $. Then we
demand the trivializing morphisms $\lambda _\alpha $ to be agree, which
means that the diagrams

\begin{equation}
\setlength{\unitlength}{.25in}\begin{picture}(4,4)
\put(-0.9,3){\makebox(1,1){\large ${E_{\alpha} \cap E_{\beta}}$}}
\put(4.5,3){\makebox(1,1){\large ${U_{\alpha} \cap U_\beta \times F}$}}
\put(3.8,0){\makebox(1,1){\large ${U_{\alpha} \cap U_\beta \times F}$}}
\put(0.7,2){{\small $\lambda_\alpha $}} \put(1,3){\vector(1,-1){1.9}}
\put(1.6,3.7){{\small $\lambda_\beta $}} \put(1.1,3.5){\vector(1,0){1.6}}
\put(3.5,2.85){\vector(0,-1){1.75}} \put(3.7,2){{\small $\Lambda_{\alpha
\beta}$}} \end{picture}  \label{p4}
\end{equation}

and

\begin{equation}
\setlength{\unitlength}{.25in}\begin{picture}(4,4)
\put(-0.9,3){\makebox(1,1){\large ${E_{\alpha} \cap E_{\beta}}$}}
\put(4.5,3){\makebox(1,1){\large ${U_{\alpha} \cap U_\beta \times F}$}}
\put(3.8,0){\makebox(1,1){\large ${U_{\alpha} \cap U_\beta \times F}$}}
\put(0.7,2){{\small $\lambda_\alpha $}} \put(1,3){\vector(1,-1){1.9}}
\put(1.6,3.7){{\small $\lambda_\beta $}} \put(1.1,3.5){\vector(1,0){1.6}}
\put(3.5,1.1){\vector(0,1){1.8}} \put(3.7,2){{\small $\Lambda_{\beta \alpha
}$}} \end{picture}  \label{p5}
\end{equation}

\noindent
should commute, where $\Lambda _{\alpha \beta }$ and $\Lambda _{\beta \alpha
}$ are maps acting along a semi-fiber $F$.

\begin{definition}
Gluing {\sl semi-transition functions }of a locally trivial semi-bundle $%
{\cal L}=\left( E,M,F,\pi \right) $ are defined by the equations
\begin{equation}
\Lambda _{\alpha \beta }\circ \lambda _\beta =\lambda _\alpha ,  \label{14}
\end{equation}
\begin{equation}
\Lambda _{\beta \alpha }\circ \lambda _\alpha =\lambda _\beta .  \label{15}
\end{equation}
\end{definition}

\begin{conjecture}
The semi-transition functions of a semi-bundle ${\cal L}$ satisfy the
following relations
\begin{equation}
\Lambda _{\alpha \beta }\circ \Lambda _{\beta \alpha }\circ \Lambda _{\alpha
\beta }=\Lambda _{\alpha \beta }  \label{16}
\end{equation}
\noindent on $U_\alpha \cap U_\beta $ overlaps and
\begin{equation}
\Lambda _{\alpha \beta }\circ \Lambda _{\beta \gamma }\circ \Lambda _{\gamma
\alpha }\circ \Lambda _{\alpha \beta }=\Lambda _{\alpha \beta },  \label{a90}
\end{equation}
\begin{equation}
\Lambda _{\beta \gamma }\circ \Lambda _{\gamma \alpha }\circ \Lambda
_{\alpha \beta }\circ \Lambda _{\beta \gamma }=\Lambda _{\beta \gamma },
\label{a91}
\end{equation}
\begin{equation}
\Lambda _{\gamma \alpha }\circ \Lambda _{\alpha \beta }\circ \Lambda _{\beta
\gamma }\circ \Lambda _{\gamma \alpha }=\Lambda _{\gamma \alpha }
\label{a92}
\end{equation}
on triple overlaps $U_\alpha \cap U_\beta \cap U_\gamma $ and
\begin{equation}
\Lambda _{\alpha \beta }\circ \Lambda _{\beta \gamma }\circ \Lambda _{\gamma
\rho }\circ \Lambda _{\rho \alpha }\circ \Lambda _{\alpha \beta }=\Lambda
_{\alpha \beta },  \label{a90w}
\end{equation}
\begin{equation}
\Lambda _{\beta \gamma }\circ \Lambda _{\gamma \rho }\circ \Lambda _{\rho
\alpha }\circ \Lambda _{\alpha \beta }\circ \Lambda _{\beta \gamma }=\Lambda
_{\beta \gamma },  \label{a91w}
\end{equation}
\begin{equation}
\Lambda _{\gamma \rho }\circ \Lambda _{\rho \alpha }\circ \Lambda _{\alpha
\beta }\circ \Lambda _{\beta \gamma }\circ \Lambda _{\gamma \rho }=\Lambda
_{\gamma \rho },  \label{a92v}
\end{equation}
\begin{equation}
\Lambda _{\rho \alpha }\circ \Lambda _{\alpha \beta }\circ \Lambda _{\beta
\gamma }\circ \Lambda _{\gamma \rho }\circ \Lambda _{\rho \alpha }=\Lambda
_{\rho \alpha }  \label{a92w}
\end{equation}
\noindent on $U_\alpha \cap U_\beta \cap U_\gamma \cap U_\rho $ .
\end{conjecture}

\begin{definition}
A semi-bundle ${\cal L}$ is called {\sl reflexive} if, in addition to (\ref
{16})-(\ref{a92w}), the semi-transition functions satisfy to the reflexivity
conditions
\begin{equation}
\Lambda _{\beta \alpha }\circ \Lambda _{\alpha \beta }\circ \Lambda _{\beta
\alpha }=\Lambda _{\beta \alpha }  \label{16a}
\end{equation}
\noindent on $U_\alpha \cap U_\beta $ overlaps and
\begin{equation}
\Lambda _{\alpha \gamma }\circ \Lambda _{\gamma \beta }\circ \Lambda _{\beta
\alpha }\circ \Lambda _{\alpha \gamma }=\Lambda _{\alpha \gamma },
\label{a9a}
\end{equation}
\begin{equation}
\Lambda _{\gamma \beta }\circ \Lambda _{\beta \alpha }\circ \Lambda _{\alpha
\gamma }\circ \Lambda _{\gamma \beta }=\Lambda _{\gamma \beta },
\label{a9a1}
\end{equation}
\begin{equation}
\Lambda _{\beta \alpha }\circ \Lambda _{\alpha \gamma }\circ \Lambda
_{\gamma \beta }\circ \Lambda _{\beta \alpha }=\Lambda _{\beta \alpha }
\label{a9a2}
\end{equation}
on triple overlaps $U_\alpha \cap U_\beta \cap U_\gamma $ and
\begin{equation}
\Lambda _{\alpha \rho }\circ \Lambda _{\rho \gamma }\circ \Lambda _{\gamma
\beta }\circ \Lambda _{\beta \alpha }\circ \Lambda _{\alpha \rho }=\Lambda
_{\alpha \rho },  \label{a9aw}
\end{equation}
\begin{equation}
\Lambda _{\rho \gamma }\circ \Lambda _{\gamma \beta }\circ \Lambda _{\beta
\alpha }\circ \Lambda _{\alpha \rho }\circ \Lambda _{\rho \gamma }=\Lambda
_{\rho \gamma },  \label{a9aw1}
\end{equation}
\begin{equation}
\Lambda _{\gamma \beta }\circ \Lambda _{\beta \alpha }\circ \Lambda _{\alpha
\rho }\circ \Lambda _{\rho \gamma }\circ \Lambda _{\gamma \beta }=\Lambda
_{\gamma \beta },  \label{a9a1w}
\end{equation}
\begin{equation}
\Lambda _{\beta \alpha }\circ \Lambda _{\alpha \rho }\circ \Lambda _{\rho
\gamma }\circ \Lambda _{\gamma \beta }\circ \Lambda _{\beta \alpha }=\Lambda
_{\beta \alpha }  \label{a9a2w}
\end{equation}

\noindent on $U_\alpha \cap U_\beta \cap U_\gamma \cap U_\rho $ .
\end{definition}

For a fixed $b\in U_\alpha \cap U_\beta $ gluing transition function $%
\Lambda _{\alpha \beta }$ describe morphisms of a semi-fiber $F$ to itself
by the condition
\begin{equation}
\Lambda _{\alpha \beta }:\left( b,f\right) \rightarrow \left( b,L_{\alpha
\beta }f\right) ,  \label{16b}
\end{equation}
where $L_{\alpha \beta }:U_\alpha \cap U_\beta \rightarrow F$ and $f\in F$.
The functions $L_{\alpha \beta }$ satisfy to the generalized cocycle
conditions similar to (\ref{16})--(\ref{a9a2w}).

\begin{remark}
The sections and transition functions of a (super) fiber bundle are
noninvertible even in the standard definitions \cite{husemoller,lang,tuy}.
But this noninvertibility has different nature as compare with our
assumptions. The transition functions are implied to be homeomorphisms and
sections should be in 1-1 correspondence with maps from the base to the
fiber \cite{porter,rab/cra1}. Our definitions (\ref{8}-\ref{9a2w}) and (\ref
{12a}-\ref{a9a2w}) extend them, allowing to include in consideration
properly noninvertible superfunctions as well.
\end{remark}

\subsection{Morphisms of semi-bundles}

Let ${\cal L}=\left( E,M,F,\pi \right) $ and ${\cal L}^{\prime }=\left(
E^{\prime },M^{\prime },F^{\prime },\pi ^{\prime }\right) $ be two
semi-bundles.

\begin{definition}
A semi-bundle morphism ${\cal L}\stackrel{f}{\rightarrow }{\cal L}^{\prime }$
consists of two morphisms $f=\left( f_E,f_M\right) $ , where $%
f_E:E\rightarrow E^{\prime }$ and $f_M:M\rightarrow M^{\prime }$ ,satisfying
$f_M\circ \pi =\pi ^{\prime }\circ f_E$ such that the diagram

\begin{equation}
\setlength{\unitlength}{.25in}%
\begin{picture}(4,4) \put(0.2,3.1){\makebox(1,1){\large ${E}$}}
\put(0.2,0){\makebox(1,1){\large
${M}$}} \put(3.1,3.05){\makebox(1,1){\large ${E^\prime}$}}
\put(3.1,0){\makebox(1,1){\large ${M^\prime}$}}
\put(0,2){{\small $\pi$}}
\put(0.7,3){\vector(0,-1){1.9}}
\put(1.8,3.8){{\small $f_E$}}
\put(1.8,1){{\small $f_M$}}
\put(1.2,3.5){\vector(1,0){1.9}}
\put(1.2,0.5){\vector(1,0){1.9}}
\put(3.5,3){\vector(0,-1){1.9}}
\put(3.7,2){{\small $\pi^\prime$}} \end{picture}
\label{p6}
\end{equation}

\noindent is commutative.
\end{definition}

Let $E_b=\left\{ x\in E\,|\,\pi \left( x\right) =b,\,b\in U\subset M\right\}
,$ then $f_E\left( E_b\right) \subset E_{f_M\left( b\right) }^{\prime }$ for
each $b$, and so the semi-fiber over $b\in M$ is carried into the semi-fiber
over $f\left( b\right) \in M^{\prime }$ by $f_E$ being a fiber morphism. If
a semi-bundle has a section, $f_E$ acts as follows $s\left( b\right)
\rightarrow s^{\prime }\left( f_M\left( b\right) \right) $.

In most applications of fiber bundles the morphism $f_M$ is identity, and $%
f_0=\left( f_E,\limfunc{id}\right) $ is called $B${\it -morphism} \cite
{husemoller}. Nevertheless, in case of semi-bundles an opposite extreme
situation can take place, when $f_M$ is a noninvertible morphism.

For each fixed $b\in M$ there exist trivializing maps $\lambda
:E_b\rightarrow U\times F$ and $\lambda ^{\prime }:E_{f_M\left( b\right)
}\rightarrow U^{\prime }\times F^{\prime }$ , $f_M\left( U\right) \subset
U^{\prime }$ which lead to a map of semi-fibers $h_b$ determined by the
commutative diagram

\begin{equation}
\setlength{\unitlength}{.25in}%
\begin{picture}(4,4) \put(0.2,3){\makebox(1,1){\large ${E_b}$}}
\put(-0.5,0){\makebox(1,1){\large
${U \times F}$}} \put(3.6,2.95){\makebox(1,1){\large ${E^\prime_{f_M \left(
b\right)}}$}} \put(4.1,0){\makebox(1,1){\large ${U^\prime \times
F^\prime}$}} \put(0,2){{\small $\lambda$}}
\put(3.7,2){{\small $\lambda^\prime$}}
\put(1.4,3.8){{\small $f_E\left( b\right)$}}
\put(1.8,0.8){{\small $h_b$}}
\put(0.7,3){\vector(0,-1){1.9}}
\put(1.2,3.5){\vector(1,0){1.9}}
\put(1.2,0.5){\vector(1,0){1.9}}
\put(3.5,3){\vector(0,-1){1.9}} \end{picture}
\label{p7}
\end{equation}

To describe a semi-bundle morphism ${\cal L}\stackrel{f}{\rightarrow }{\cal L%
}^{\prime }$ locally we choose open coverings $M=\cup U_\alpha $ and $%
M^{\prime }=\cup U_{\alpha ^{\prime }}^{\prime }$ together with
trivializations $\lambda _\alpha $ and $\lambda _{\alpha ^{\prime }}^{\prime
}$ (see (\ref{p3})). Then the connection between semi-transition functions $%
\Lambda _{\alpha \beta }$ and $\Lambda _{\alpha ^{\prime }\beta ^{\prime
}}^{\prime }$ (\ref{14})--(\ref{15}) of two semi-bundles ${\cal L}$ and $%
{\cal L}^{\prime }$ can be found from the commutative diagram

\begin{equation}
\setlength{\unitlength}{.25in}
\begin{picture}(4,4) \put(-0.8,3){\makebox(1,1){\large ${U_{\alpha \beta}
\times F}$}} \put(-1.05,0){\makebox(1,1){\large ${U^\prime_{\alpha^\prime
\beta^\prime}
\times F^\prime}$}} \put(4.2,3){\makebox(1,1){\large ${U_{\alpha \beta}
\times F}$}}
\put(4.4,0){\makebox(1,1){\large ${U^\prime_{\alpha^\prime \beta^\prime}
\times F^\prime}$}} \put(-0.1,2){{\small $h_\alpha$}}
\put(3.7,2){{\small $h_\beta$}}
\put(1.8,3.8){{\small $\Lambda_{\alpha \beta}$}}
\put(1.8,1){{\small $\Lambda^\prime_{\alpha^\prime \beta^\prime}$}}
\put(0.7,3){\vector(0,-1){1.9}}
\put(3.5,3){\vector(0,-1){1.9}} \put(1.2,3.5){\vector(1,0){1.9}}
\put(1.2,0.6){\vector(1,0){1.9}} \end{picture}
\label{p8}
\end{equation}

\noindent where morphisms $h_\alpha $ are defined by the diagram %

\begin{equation}
\setlength{\unitlength}{.25in}%
\begin{picture}(4,4) \put(0.25,3){\makebox(1,1){\large ${E}$}}
\put(-0.7,0){\makebox(1,1){\large
${U_{\alpha} \times F}$}} \put(3.1,3){\makebox(1,1){\large ${E^\prime}$}}
\put(4.2,0){\makebox(1,1){\large ${U^\prime_{\alpha^\prime} \times
F^\prime}$}} \put(-0.1,2){{\small $\lambda_\alpha$}}
\put(3.7,2){{\small $\lambda^\prime_{\alpha^\prime}$}}
\put(1.8,3.8){{\small $f_E$}}
\put(1.8,0.8){{\small $h_\alpha$}}
\put(0.7,3){\vector(0,-1){1.9}}
\put(1.2,3.5){\vector(1,0){1.9}}
\put(1.2,0.5){\vector(1,0){1.9}} \put(3.5,3){\vector(0,-1){1.9}}
\end{picture}
\label{p9}
\end{equation}

From (\ref{p8}) we have the relation between semi-transition functions
\begin{equation}
h_\alpha \circ \Lambda _{\alpha \beta }=\Lambda _{\alpha ^{\prime }\beta
^{\prime }}^{\prime }\circ h_\beta  \label{17}
\end{equation}
which holds valid for noninvertible $h_\alpha $ as well, while in the
invertible case \cite{husemoller,lang} the equation (\ref{17}) is solved
with respect to $\Lambda _{\alpha ^{\prime }\beta ^{\prime }}^{\prime }$ ,
as follows $\Lambda _{\alpha ^{\prime }\beta ^{\prime }}^{\prime }=h_\alpha
\circ \Lambda _{\alpha \beta }\circ h_\beta ^{-1}$ (it can be considered as
an equivalence of cocycles). However, in general (\ref{17}) is a system of
superequations which should be solved by the standard \cite{berezin} or
extended \cite{ber/rab} methods of superanalysis.

Let $M$ admits two trivializing coverings $\left\{ U_\alpha ,\lambda _\alpha
\right\} $ and $\left\{ U_{\alpha ^{\prime }}^{\prime },\lambda _{\alpha
^{\prime }}^{\prime }\right\} $. In general they are not connected and
semi-transition functions $\Lambda _{\alpha \beta }$ and $\Lambda _{\alpha
^{\prime }\beta ^{\prime }}^{\prime }$ are independent. However, if $M$ is
the base superspace for two semi-bundles ${\cal L}$ and ${\cal L}^{\prime }$
which are connected by a $B$-morphism ${\cal L}\stackrel{f_0}{\rightarrow }%
{\cal L}^{\prime }$ , then $\Lambda _{\alpha \beta }$ and $\Lambda _{\alpha
^{\prime }\beta ^{\prime }}^{\prime }$ should agree properly.

\begin{proposition}
The semi-transition functions $\Lambda _{\alpha \beta }$ and $\Lambda
_{\alpha ^{\prime }\beta ^{\prime }}^{\prime }$ agree if there exist
additional maps $\tilde{\Lambda}_{\alpha ^{\prime }\beta }:U_{\alpha
^{\prime }}^{\prime }\cap U_\beta $ and $\tilde{\Lambda}_{\alpha \beta
^{\prime }}:U_\alpha \cap U_{\beta ^{\prime }}^{\prime }$ connected between
themselves by the relations
\begin{equation}
\tilde{\Lambda}_{\alpha ^{\prime }\beta }\circ \tilde{\Lambda}_{\beta \alpha
^{\prime }}\circ \tilde{\Lambda}_{\alpha ^{\prime }\beta }=\tilde{\Lambda}%
_{\alpha ^{\prime }\beta }  \label{16c}
\end{equation}
on $U_{\alpha ^{\prime }}^{\prime }\cap U_\beta $ and
\begin{equation}
\tilde{\Lambda}_{\alpha \beta ^{\prime }}\circ \tilde{\Lambda}_{\beta
^{\prime }\alpha }\circ \tilde{\Lambda}_{\alpha \beta ^{\prime }}=\tilde{%
\Lambda}_{\alpha \beta ^{\prime }}  \label{16d}
\end{equation}
on $U_\alpha \cap U_{\beta ^{\prime }}^{\prime }$ overlaps.

Then the agreement conditions for $\Lambda _{\alpha \beta }$ and $\Lambda
_{\alpha ^{\prime }\beta ^{\prime }}^{\prime }$ are
\begin{equation}
\tilde{\Lambda}_{\alpha ^{\prime }\beta }\circ \Lambda _{\beta \gamma }\circ
\tilde{\Lambda}_{\gamma \alpha ^{\prime }}\circ \tilde{\Lambda}_{\alpha
^{\prime }\beta }=\tilde{\Lambda}_{\alpha ^{\prime }\beta },  \label{18a}
\end{equation}
\begin{equation}
\Lambda _{\beta \gamma }\circ \tilde{\Lambda}_{\gamma \alpha ^{\prime
}}\circ \tilde{\Lambda}_{\alpha ^{\prime }\beta }\circ \Lambda _{\beta
\gamma }=\Lambda _{\beta \gamma },  \label{18b}
\end{equation}
\begin{equation}
\tilde{\Lambda}_{\gamma \alpha ^{\prime }}\circ \tilde{\Lambda}_{\alpha
^{\prime }\beta }\circ \Lambda _{\beta \gamma }\circ \tilde{\Lambda}_{\gamma
\alpha ^{\prime }}=\tilde{\Lambda}_{\gamma \alpha ^{\prime }}  \label{18c}
\end{equation}
on triple overlaps $U_{\alpha ^{\prime }}^{\prime }\cap U_\beta \cap
U_\gamma $ and
\begin{equation}
\Lambda _{\alpha ^{\prime }\beta ^{\prime }}^{\prime }\circ \tilde{\Lambda}%
_{\beta ^{\prime }\gamma }\circ \tilde{\Lambda}_{\gamma \alpha ^{\prime
}}\circ \Lambda _{\alpha ^{\prime }\beta ^{\prime }}^{\prime }=\Lambda
_{\alpha ^{\prime }\beta ^{\prime }}^{\prime },  \label{19a}
\end{equation}
\begin{equation}
\tilde{\Lambda}_{\beta ^{\prime }\gamma }\circ \tilde{\Lambda}_{\gamma
\alpha ^{\prime }}\circ \Lambda _{\alpha ^{\prime }\beta ^{\prime }}^{\prime
}\circ \tilde{\Lambda}_{\beta ^{\prime }\gamma }=\tilde{\Lambda}_{\beta
^{\prime }\gamma },  \label{19b}
\end{equation}
\begin{equation}
\tilde{\Lambda}_{\gamma \alpha ^{\prime }}\circ \Lambda _{\alpha ^{\prime
}\beta ^{\prime }}^{\prime }\circ \tilde{\Lambda}_{\beta ^{\prime }\gamma
}\circ \tilde{\Lambda}_{\gamma \alpha ^{\prime }}=\tilde{\Lambda}_{\gamma
\alpha ^{\prime }}  \label{19c}
\end{equation}
on $U_{\alpha ^{\prime }}^{\prime }\cap U_{\beta ^{\prime }}^{\prime }\cap
U_\gamma $ overlaps. Then
\begin{equation}
\tilde{\Lambda}_{\alpha ^{\prime }\beta }\circ \Lambda _{\beta \gamma }\circ
\Lambda _{\gamma \rho }\circ \tilde{\Lambda}_{\rho \alpha ^{\prime }}\circ
\tilde{\Lambda}_{\alpha ^{\prime }\beta }=\tilde{\Lambda}_{\alpha ^{\prime
}\beta },  \label{ka90w}
\end{equation}
\begin{equation}
\Lambda _{\beta \gamma }\circ \Lambda _{\gamma \rho }\circ \tilde{\Lambda}%
_{\rho \alpha ^{\prime }}\circ \tilde{\Lambda}_{\alpha ^{\prime }\beta
}\circ \Lambda _{\beta \gamma }=\Lambda _{\beta \gamma },  \label{ka91w}
\end{equation}
\begin{equation}
\Lambda _{\gamma \rho }\circ \tilde{\Lambda}_{\rho \alpha ^{\prime }}\circ
\tilde{\Lambda}_{\alpha ^{\prime }\beta }\circ \Lambda _{\beta \gamma }\circ
\Lambda _{\gamma \rho }=\Lambda _{\gamma \rho },  \label{ka92v}
\end{equation}
\begin{equation}
\tilde{\Lambda}_{\rho \alpha ^{\prime }}\circ \tilde{\Lambda}_{\alpha
^{\prime }\beta }\circ \Lambda _{\beta \gamma }\circ \Lambda _{\gamma \rho
}\circ \tilde{\Lambda}_{\rho \alpha ^{\prime }}=\tilde{\Lambda}_{\rho \alpha
^{\prime }}  \label{ka92w}
\end{equation}
on $U_{\alpha ^{\prime }}^{\prime }\cap U_\beta \cap U_\gamma \cap U_\rho $
and
\begin{equation}
\Lambda _{\alpha ^{\prime }\beta ^{\prime }}^{\prime }\circ \tilde{\Lambda}%
_{\beta ^{\prime }\gamma }\circ \Lambda _{\gamma \rho }\circ \tilde{\Lambda}%
_{\rho \alpha ^{\prime }}\circ \tilde{\Lambda}_{\alpha \beta ^{\prime
}}=\Lambda _{\alpha ^{\prime }\beta ^{\prime }}^{\prime },  \label{ma90w}
\end{equation}
\begin{equation}
\tilde{\Lambda}_{\beta ^{\prime }\gamma }\circ \Lambda _{\gamma \rho }\circ
\tilde{\Lambda}_{\rho \alpha ^{\prime }}\circ \Lambda _{\alpha ^{\prime
}\beta ^{\prime }}^{\prime }\circ \tilde{\Lambda}_{\beta ^{\prime }\gamma }=%
\tilde{\Lambda}_{\beta ^{\prime }\gamma },  \label{ma91w}
\end{equation}
\begin{equation}
\Lambda _{\gamma \rho }\circ \tilde{\Lambda}_{\rho \alpha ^{\prime }}\circ
\Lambda _{\alpha ^{\prime }\beta ^{\prime }}^{\prime }\circ \tilde{\Lambda}%
_{\beta ^{\prime }\gamma }\circ \Lambda _{\gamma \rho }=\Lambda _{\gamma
\rho },  \label{ma92v}
\end{equation}
\begin{equation}
\tilde{\Lambda}_{\rho \alpha ^{\prime }}\circ \Lambda _{\alpha ^{\prime
}\beta ^{\prime }}^{\prime }\circ \tilde{\Lambda}_{\beta ^{\prime }\gamma
}\circ \Lambda _{\gamma \rho }\circ \tilde{\Lambda}_{\rho \alpha ^{\prime }}=%
\tilde{\Lambda}_{\rho \alpha ^{\prime }}  \label{ma92w}
\end{equation}
on $U_{\alpha ^{\prime }}^{\prime }\cap U_{\beta ^{\prime }}^{\prime }\cap
U_\gamma \cap U_\rho $ and
\begin{equation}
\Lambda _{\alpha ^{\prime }\beta ^{\prime }}^{\prime }\circ \Lambda _{\beta
^{\prime }\gamma ^{\prime }}^{\prime }\circ \tilde{\Lambda}_{\gamma ^{\prime
}\rho }\circ \tilde{\Lambda}_{\rho \alpha ^{\prime }}\circ \Lambda _{\alpha
^{\prime }\beta ^{\prime }}^{\prime }=\Lambda _{\alpha ^{\prime }\beta
^{\prime }}^{\prime },  \label{na90w}
\end{equation}
\begin{equation}
\Lambda _{\beta ^{\prime }\gamma ^{\prime }}^{\prime }\circ \tilde{\Lambda}%
_{\gamma ^{\prime }\rho }\circ \tilde{\Lambda}_{\rho \alpha ^{\prime }}\circ
\Lambda _{\alpha ^{\prime }\beta ^{\prime }}^{\prime }\circ \tilde{\Lambda}%
_{\beta ^{\prime }\gamma }=\tilde{\Lambda}_{\beta ^{\prime }\gamma },
\label{na91w}
\end{equation}
\begin{equation}
\tilde{\Lambda}_{\gamma ^{\prime }\rho }\circ \tilde{\Lambda}_{\rho \alpha
^{\prime }}\circ \Lambda _{\alpha ^{\prime }\beta ^{\prime }}^{\prime }\circ
\Lambda _{\beta ^{\prime }\gamma ^{\prime }}^{\prime }\circ \tilde{\Lambda}%
_{\gamma ^{\prime }\rho }=\tilde{\Lambda}_{\gamma ^{\prime }\rho },
\label{na92v}
\end{equation}
\begin{equation}
\tilde{\Lambda}_{\rho \alpha ^{\prime }}\circ \Lambda _{\alpha ^{\prime
}\beta ^{\prime }}^{\prime }\circ \Lambda _{\beta ^{\prime }\gamma ^{\prime
}}^{\prime }\circ \tilde{\Lambda}_{\gamma ^{\prime }\rho }\circ \tilde{%
\Lambda}_{\rho \alpha ^{\prime }}=\tilde{\Lambda}_{\rho \alpha ^{\prime }}
\label{na92w}
\end{equation}
on $U_{\alpha ^{\prime }}^{\prime }\cap U_{\beta ^{\prime }}^{\prime }\cap
U_{\gamma ^{\prime }}^{\prime }\cap U_\rho $.
\end{proposition}

\begin{proof}
Construct a sum of trivializing coverings $\left\{ U_\alpha ,\lambda _\alpha
\right\} $ and $\left\{ U_{\alpha ^{\prime }}^{\prime },\lambda _{\alpha
^{\prime }}^{\prime }\right\} $ and then use (\ref{16})--(\ref{a92w}).
\end{proof}

\begin{proposition}
The semi-transition functions $\Lambda _{\alpha \beta }$ and $\Lambda
_{\alpha ^{\prime }\beta ^{\prime }}^{\prime }$ reflexively agree if there
exist additional reflexive maps $\tilde{\Lambda}_{\alpha ^{\prime }\beta
}:U_{\alpha ^{\prime }}^{\prime }\cap U_\beta $ and $\tilde{\Lambda}_{\alpha
\beta ^{\prime }}:U_\alpha \cap U_{\beta ^{\prime }}^{\prime }$ connected
between themselves (in addition to (\ref{16c})--(\ref{16d})) by the
reflexive relations
\begin{equation}
\tilde{\Lambda}_{\beta \alpha ^{\prime }}\circ \tilde{\Lambda}_{\alpha
^{\prime }\beta }\circ \tilde{\Lambda}_{\beta \alpha ^{\prime }}=\tilde{%
\Lambda}_{\beta \alpha ^{\prime }}  \label{r16c}
\end{equation}
on $U_{\alpha ^{\prime }}^{\prime }\cap U_\beta $ and
\begin{equation}
\tilde{\Lambda}_{\beta ^{\prime }\alpha }\circ \tilde{\Lambda}_{\alpha \beta
^{\prime }}\circ \tilde{\Lambda}_{\beta ^{\prime }\alpha }=\tilde{\Lambda}%
_{\beta ^{\prime }\alpha }  \label{r16d}
\end{equation}
on $U_\alpha \cap U_{\beta ^{\prime }}^{\prime }$ overlaps. The reflexive
semi-transition functions $\Lambda _{\alpha \beta }$ and $\Lambda _{\alpha
^{\prime }\beta ^{\prime }}^{\prime }$ should satisfy (in addition to (\ref
{18a})--(\ref{na92w})) the following reflexivity agreement relations
\begin{equation}
\tilde{\Lambda}_{\alpha ^{\prime }\gamma }\circ \Lambda _{\gamma \beta
}\circ \tilde{\Lambda}_{\beta \alpha ^{\prime }}\circ \tilde{\Lambda}%
_{\alpha ^{\prime }\gamma }=\tilde{\Lambda}_{\alpha ^{\prime }\gamma },
\label{21a}
\end{equation}
\begin{equation}
\Lambda _{\gamma \beta }\circ \tilde{\Lambda}_{\beta \alpha ^{\prime }}\circ
\tilde{\Lambda}_{\alpha ^{\prime }\gamma }\circ \Lambda _{\gamma \beta
}=\Lambda _{\gamma \beta },  \label{21b}
\end{equation}
\begin{equation}
\tilde{\Lambda}_{\beta \alpha ^{\prime }}\circ \tilde{\Lambda}_{\alpha
^{\prime }\gamma }\circ \Lambda _{\gamma \beta }\circ \tilde{\Lambda}_{\beta
\alpha ^{\prime }}=\tilde{\Lambda}_{\beta \alpha ^{\prime }}  \label{21c}
\end{equation}
on $U_{\alpha ^{\prime }}^{\prime }\cap U_\beta \cap U_\gamma $ and
\begin{equation}
\tilde{\Lambda}_{\alpha ^{\prime }\gamma }\circ \tilde{\Lambda}_{\gamma
\beta ^{\prime }}\circ \Lambda _{\beta ^{\prime }\alpha ^{\prime }}^{\prime
}\circ \tilde{\Lambda}_{\alpha ^{\prime }\gamma }=\tilde{\Lambda}_{\alpha
^{\prime }\gamma },  \label{22a}
\end{equation}
\begin{equation}
\tilde{\Lambda}_{\gamma \beta ^{\prime }}\circ \Lambda _{\beta ^{\prime
}\alpha ^{\prime }}^{\prime }\circ \tilde{\Lambda}_{\alpha ^{\prime }\gamma
}\circ \tilde{\Lambda}_{\gamma \beta ^{\prime }}=\tilde{\Lambda}_{\gamma
\beta ^{\prime }},  \label{22b}
\end{equation}
\begin{equation}
\Lambda _{\beta ^{\prime }\alpha ^{\prime }}^{\prime }\circ \tilde{\Lambda}%
_{\alpha ^{\prime }\gamma }\circ \Lambda _{\gamma \beta ^{\prime }}^{\prime
}\circ \Lambda _{\beta ^{\prime }\alpha ^{\prime }}^{\prime }=\Lambda
_{\beta ^{\prime }\alpha ^{\prime }}^{\prime }  \label{22c}
\end{equation}
on $U_{\alpha ^{\prime }}^{\prime }\cap U_{\beta ^{\prime }}^{\prime }\cap
U_\gamma $ overlaps. Then
\begin{equation}
\tilde{\Lambda}_{\alpha ^{\prime }\rho }\circ \Lambda _{\rho \gamma }\circ
\Lambda _{\gamma \beta }\circ \tilde{\Lambda}_{\beta \alpha ^{\prime }}\circ
\tilde{\Lambda}_{\alpha ^{\prime }\rho }=\tilde{\Lambda}_{\alpha ^{\prime
}\rho },  \label{ta9aw}
\end{equation}
\begin{equation}
\Lambda _{\rho \gamma }\circ \Lambda _{\gamma \beta }\circ \tilde{\Lambda}%
_{\beta \alpha ^{\prime }}\circ \tilde{\Lambda}_{\alpha ^{\prime }\rho
}\circ \Lambda _{\rho \gamma }=\Lambda _{\rho \gamma },  \label{ta9aw1}
\end{equation}
\begin{equation}
\Lambda _{\gamma \beta }\circ \tilde{\Lambda}_{\beta \alpha ^{\prime }}\circ
\tilde{\Lambda}_{\alpha ^{\prime }\rho }\circ \Lambda _{\rho \gamma }\circ
\Lambda _{\gamma \beta }=\Lambda _{\gamma \beta },  \label{ta9a1w}
\end{equation}
\begin{equation}
\tilde{\Lambda}_{\beta \alpha ^{\prime }}\circ \tilde{\Lambda}_{\alpha
^{\prime }\rho }\circ \Lambda _{\rho \gamma }\circ \Lambda _{\gamma \beta
}\circ \tilde{\Lambda}_{\beta \alpha ^{\prime }}=\tilde{\Lambda}_{\beta
\alpha ^{\prime }}  \label{ta9a2w}
\end{equation}

\noindent on $U_{\alpha ^{\prime }}^{\prime }\cap U_\beta \cap U_\gamma \cap
U_\rho $ and
\begin{equation}
\tilde{\Lambda}_{\alpha ^{\prime }\rho }\circ \Lambda _{\rho \gamma }\circ
\tilde{\Lambda}_{\gamma \beta ^{\prime }}\circ \Lambda _{\beta ^{\prime
}\alpha ^{\prime }}^{\prime }\circ \tilde{\Lambda}_{\alpha ^{\prime }\rho }=%
\tilde{\Lambda}_{\alpha ^{\prime }\rho },  \label{ua9aw}
\end{equation}
\begin{equation}
\Lambda _{\rho \gamma }\circ \tilde{\Lambda}_{\gamma \beta ^{\prime }}\circ
\Lambda _{\beta ^{\prime }\alpha ^{\prime }}^{\prime }\circ \tilde{\Lambda}%
_{\alpha ^{\prime }\rho }\circ \Lambda _{\rho \gamma }=\Lambda _{\rho \gamma
},  \label{ua9aw1}
\end{equation}
\begin{equation}
\tilde{\Lambda}_{\gamma \beta ^{\prime }}\circ \Lambda _{\beta ^{\prime
}\alpha ^{\prime }}^{\prime }\circ \tilde{\Lambda}_{\alpha ^{\prime }\rho
}\circ \Lambda _{\rho \gamma }\circ \tilde{\Lambda}_{\gamma \beta ^{\prime
}}=\tilde{\Lambda}_{\gamma \beta ^{\prime }},  \label{ua9a1w}
\end{equation}
\begin{equation}
\Lambda _{\beta ^{\prime }\alpha ^{\prime }}^{\prime }\circ \tilde{\Lambda}%
_{\alpha ^{\prime }\rho }\circ \Lambda _{\rho \gamma }\circ \tilde{\Lambda}%
_{\gamma \beta ^{\prime }}\circ \Lambda _{\beta ^{\prime }\alpha ^{\prime
}}^{\prime }=\Lambda _{\beta ^{\prime }\alpha ^{\prime }}^{\prime }
\label{ua9a2w}
\end{equation}

\noindent on $U_{\alpha ^{\prime }}^{\prime }\cap U_{\beta ^{\prime
}}^{\prime }\cap U_\gamma \cap U_\rho $ and
\begin{equation}
\tilde{\Lambda}_{\alpha ^{\prime }\rho }\circ \tilde{\Lambda}_{\rho \gamma
^{\prime }}\circ \Lambda _{\gamma ^{\prime }\beta ^{\prime }}^{\prime }\circ
\Lambda _{\beta ^{\prime }\alpha ^{\prime }}^{\prime }\circ \tilde{\Lambda}%
_{\alpha ^{\prime }\rho }=\tilde{\Lambda}_{\alpha ^{\prime }\rho },
\label{ya9aw}
\end{equation}
\begin{equation}
\tilde{\Lambda}_{\rho \gamma ^{\prime }}\circ \Lambda _{\gamma ^{\prime
}\beta ^{\prime }}^{\prime }\circ \Lambda _{\beta ^{\prime }\alpha ^{\prime
}}^{\prime }\circ \tilde{\Lambda}_{\alpha ^{\prime }\rho }\circ \tilde{%
\Lambda}_{\rho \gamma ^{\prime }}=\tilde{\Lambda}_{\rho \gamma ^{\prime }},
\label{ya9aw1}
\end{equation}
\begin{equation}
\Lambda _{\gamma ^{\prime }\beta ^{\prime }}^{\prime }\circ \Lambda _{\beta
^{\prime }\alpha ^{\prime }}^{\prime }\circ \tilde{\Lambda}_{\alpha ^{\prime
}\rho }\circ \tilde{\Lambda}_{\rho \gamma ^{\prime }}\circ \Lambda _{\gamma
^{\prime }\beta ^{\prime }}^{\prime }=\Lambda _{\gamma ^{\prime }\beta
^{\prime }}^{\prime },  \label{ya9a1w}
\end{equation}
\begin{equation}
\Lambda _{\beta ^{\prime }\alpha ^{\prime }}^{\prime }\circ \tilde{\Lambda}%
_{\alpha ^{\prime }\rho }\circ \tilde{\Lambda}_{\rho \gamma ^{\prime }}\circ
\Lambda _{\gamma ^{\prime }\beta ^{\prime }}^{\prime }\circ \Lambda _{\beta
^{\prime }\alpha ^{\prime }}^{\prime }=\Lambda _{\beta ^{\prime }\alpha
^{\prime }}^{\prime }  \label{ya9a2w}
\end{equation}

\noindent on $U_{\alpha ^{\prime }}^{\prime }\cap U_{\beta ^{\prime
}}^{\prime }\cap U_{\gamma ^{\prime }}^{\prime }\cap U_\rho $.
\end{proposition}

Analogously we can define and study a principal and associated semi-bundles
with a structure semigroup, but this is a subject of a separate paper which
will appear elsewhere.

\section{Noninvertibility and semi-homotopies}

Here we briefly dwell on some possibilities to extend directly the notion of
homotopy to noninvertible continuous mappings.

A {\it homotopy} \cite{fom/fuc/gut,maclane} is a continuous mapping between
two maps $f:X\rightarrow Y$ and $g:X\rightarrow Y$ in the space $C\left(
X,Y\right) $ of maps $X\rightarrow Y$ such that $\gamma _{t=0}\left(
x\right) =f\left( x\right) ,\;\gamma _{t=1}\left( x\right) =g\left( x\right)
$. Such maps are called {\it homotopic}. In other words \cite{switzer} a
homotopy from $X$ to $Y$ is a continuous function $\Gamma :X\times
I\rightarrow Y$ where $I=\left[ 0,1\right] $ is a unit interval. For a fixed
$t\in I$ one has {\it stages} $\gamma _t:X\rightarrow Y$ defined by $\gamma
_t\left( x\right) =\Gamma \left( x,t\right) $. The relation of homotopy
divides $C\left( X,Y\right) $ into a set of equivalent classes $\pi \left(
X,Y\right) $ called {\it homotopy classes} which are a set of connected
components of $C\left( X,Y\right) $. Therefore, $\pi \left( *,Y\right) $ ($*$
is a point) the homotopy classes correspond to connected components of $Y$.
If $C\left( X,Y\right) $ is connected, then the homotopy between $f\left(
x\right) $ and $g\left( x\right) $ can be chosen as their average, i.e.
\begin{equation}
\gamma _t\left( x\right) =tf\left( x\right) +\left( 1-t\right) g\left(
x\right) .  \label{23c}
\end{equation}

Two maps $f$ and $g$ are {\it homotopically equivalent} if $f\circ g$ and $%
g\circ f$ are homotopic to the identity.

Now let $X$ and $Y$ are supermanifolds in some of the definitions \cite
{bry2,cra/rab,rog1} or semi-supermanifold in our above definition. Then
there exist a possibility to extend the notion of homotopy. The idea is in
extending the definition of the parameter $t$. In the standard case the unit
interval $I=\left[ 0,1\right] $ was taken for simplicity, because any two
intervals on a real axis are homeomorphic, and so they are topologically
equal.

In supercase the situation is totally different. Let $X$ and $Y$ are defined
over $\Lambda $ , a commutative ${\Bbb Z}_2$-graded superalgebra admitting a
decomposition into direct sum $\Lambda =\Lambda _0\oplus \Lambda _1$ of the
even $\Lambda _0$ and odd $\Lambda _1$ parts and into the direct sum $%
\Lambda ={\Bbb B\oplus S}$ of the body ${\Bbb B}$ and soul ${\Bbb S}$ (see
\cite{rab/cra2,rog1} for details). The body map $\varepsilon :\Lambda
\rightarrow {\Bbb B}$ can be viewed as discarding all nilpotent superalgebra
generators, which gives a number part. So we have three topologically
disjoint cases:

\begin{enumerate}
\item  The parameter $t\in \Lambda _0$ is even and has a body, i.e. $%
\varepsilon \left( t\right) \neq 0$.

\item  The parameter $t\in \Lambda _0$ is even and has no body, i.e. $%
\varepsilon \left( t\right) =0$.

\item  The parameter $\tau \in \Lambda _1$ is odd (any odd element has no
body).
\end{enumerate}

The first choice can be reduced to the standard case by means of a
corresponding homeomorphism, and such $t$ can always be considered in the
unit interval $I=\left[ 0,1\right] $ . However, the following two
possibilities are topologically disjoint from the first one and between
themselves.

\begin{definition}
An {\sl even semi-homotopy} between two supermaps $f:X\rightarrow Y$ and $%
g:X\rightarrow Y$ is a noninvertible (in general) mapping $X\rightarrow Y$
depending on a nilpotent bodyless even parameter $t\in \Lambda _0$ and two
bodyless even constants $a,b\in \Lambda _0$ such that
\begin{equation}
\begin{array}{c}
\Delta I^{ab}\gamma _{t=a}^{even}=\Delta I^{ab}f\left( x\right) , \\
\Delta I^{ab}\gamma _{t=b}^{even}=\Delta I^{ab}g\left( x\right) ,
\end{array}
\label{23a}
\end{equation}

\noindent where
\begin{equation}
\begin{array}{c}
\gamma _t^{even}\left( x\right) =\Gamma ^{even}\left( x,t\right) ,\;\Gamma
^{even}:X\times I^{ab}\rightarrow Y, \\
I^{ab}=\left[ a,b\right] ,\;\Delta I^{ab}=b-a.
\end{array}
\label{23b}
\end{equation}
\end{definition}

\begin{definition}
An {\sl odd semi-homotopy} between two supermaps $f:X\rightarrow Y$ and $%
g:X\rightarrow Y$ is a noninvertible (in general) mapping $X\rightarrow Y$
depending on a nilpotent odd parameter $\tau \in \Lambda _1$ and two odd
constants $\mu ,\nu \in \Lambda _1$ such that
\begin{equation}
\begin{array}{c}
\Delta {\cal I}^{\alpha \beta }\gamma _{\tau =\alpha }^{odd}=\Delta {\cal I}%
^{\alpha \beta }f\left( x\right) , \\
\Delta {\cal I}^{\alpha \beta }\gamma _{\tau =\beta }^{odd}=\Delta {\cal I}%
^{\alpha \beta }g\left( x\right) .
\end{array}
\label{24a}
\end{equation}
\begin{equation}
\begin{array}{c}
\gamma _\tau ^{odd}\left( x\right) =\Gamma ^{odd}\left( x,\tau \right)
,\Gamma ^{odd}:X\times {\cal I}^{\alpha \beta }\rightarrow Y, \\
{\cal I}^{\alpha \beta }=\left[ \alpha ,\beta \right] ,\;\Delta {\cal I}%
^{\alpha \beta }=\beta -\alpha .
\end{array}
\label{24b}
\end{equation}
\end{definition}

\begin{remark}
In (\ref{23b}) and (\ref{24b}) $I^{ab}$ and ${\cal I}^{\alpha \beta }$ are
not intervals in any sense, because among bodyless variables there is no
possibility to establish an order relation \cite{bry2,cat/rei/teo,rab/cra1},
and so $\Delta I^{ab}$ and $\Delta {\cal I}^{\alpha \beta }$ are only
notations.
\end{remark}

Nevertheless, we can give an example of an analog of the average (\ref{23c})
for an odd semi-homotopy
\begin{equation}
\left( \beta -\alpha \right) \gamma _\tau ^{odd}\left( x\right) =\left(
\beta -\tau \right) f\left( x\right) +\left( \tau -\alpha \right) g\left(
x\right)  \label{25}
\end{equation}
which can satisfy some supersmooth conditions.

\begin{remark}
In (\ref{23a}) and (\ref{24a}) it is not possible to cancel the left and
right hand parts by $I^{ab}$ and ${\cal I}^{\alpha \beta }$ correspondingly,
because the solutions for semi-homotopies $\gamma _t^{even}$ and $\gamma
_\tau ^{odd}$ are viewed as equivalence relations. This is clearly seen from
(\ref{25}) where the division by $\left( \beta -\alpha \right) $ is
impossible, nevertheless a solution for $\gamma _\tau ^{odd}\left( x\right) $
exists.
\end{remark}

The most important property of semi-homotopies is their possible
noninvertibility which follows from the nilpotency of $t$ and $\tau $ and
the definitions (\ref{23a}) and (\ref{24a}). Therefore, $Y$ cannot be a
supermanifold, it can be a semi-supermanifold only.

It can be assumed that semi-homotopies play the same or similar role in the
study of continuous properties and classification of semi-supermanifolds, as
ordinary homotopies play for ordinary manifolds. So that it is worthwhile to
study their properties thoroughly and in more detail, which will be done
elsewhere.

\section{Acknowledgments}

The author is grateful to P. Van Nieuwenhuizen, J. M. Rabin and J. Stasheff
for fruitful conversations and valuable remarks. Also the discussions with
C. Day, M. Duff, J. Kupsch, M. Markl, J. McCleary, R. Mohapatra, W.
R\"{u}hl, R. Umble, A. Voronov, G. Zuckerman and B. Zwiebach are greatly
acknowledged.

\end{document}